\magnification=\magstep1
\baselineskip=14pt
\font\tiny=cmr8
\font\Bbb=msbm10
\input epsf

\baselineskip=14pt

\def\bdy{\partial}
\def\:{\colon}
\def\P{{\cal P}}

\def\R{\hbox{\Bbb R}}
 
\def\qed{\hfill \rlap{$\sqcup$}$\sqcap$}
\def\cond#1{\par\noindent\rlap{#1}\hskip35pt\hangindent=35pt\hangafter=1}

\centerline{\bf Pants Decompositions of Surfaces} 
\vskip .1truein
\centerline{Allen Hatcher}
\vskip .4truein

In studying the geometry and topology of surfaces it often happens that one considers a 
collection of disjointly embedded circles in a compact orientable surface $\Sigma$ which
decompose
$\Sigma$ into pairs of pants ---  surfaces of genus zero with three boundary circles. If
 $\Sigma$ is not itself a pair of pants, then there
are infinitely many different isotopy classes of pants decompositions of $\Sigma$. It was
observed in [HT] that any two isotopy classes of pants decompositions can be joined by a
finite sequence of ``elementary moves" in which only one circle changes at a time. In the
present paper we apply the techniques of [HT] to study the relations which hold among
such sequences of elementary moves. The main result is that there are five basic types of
relations from which all others follow. Namely, we construct a two-dimensional cell
complex $\P(\Sigma)$ whose vertices are the isotopy classes of pants decompositions of
$\Sigma$, whose edges are the elementary moves, and whose 2-cells are attached via the
basic relations. Then we prove that $\P(\Sigma)$ is simply-connected.

Now let us give the precise definitions. Let $\Sigma$ be a connected compact orientable surface. We say $\Sigma$ has  type $ (g,n) $ if it has genus $ g $ and $ n $ boundary components. By a  {\it pants decomposition\/} of $ \Sigma $ we mean a finite collection $P $ of disjoint smoothly embedded circles cutting $ \Sigma $ into pieces  which are surfaces of type $ (0,3) $. We also call $ P $ a {\it maximal cut system}. The number of curves in a maximal cut system is $ 3g - 3 + n $, and the number of complementary components is $ 2g - 2 + n = |\chi(\Sigma)| $, assuming that $ \Sigma $ has at least one pants decomposition. 

Let $P $ be a pants decomposition, and  suppose that one of the circles $ \beta $ of $
P $ is such that deleting $ \beta $ from $ P $ produces a complementary component of
type $ (1,1) $.  This is equivalent to saying there is a circle $ \gamma $ in $\Sigma$
which intersects  $ \beta $ in one point transversely and is disjoint from all the other
circles  in $ P $. In this case, replacing $ \beta $ by $ \gamma $ in $ P $ produces 
a new pants decomposition $ P'$.  We call this replacement a  {\it simple move\/}, or
S-move.
\vskip.2cm
\centerline{\epsfbox{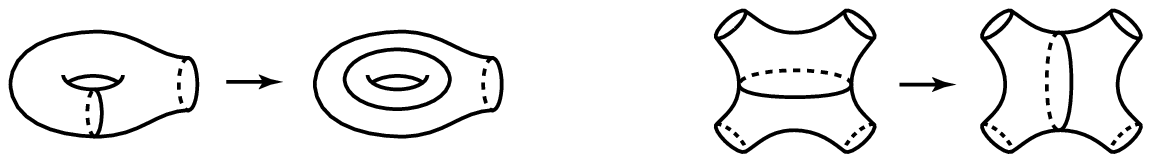}}
\vskip-.3cm
\centerline {{\bf Figure 1}: an S-move and an A-move}
 \vskip .2cm\noindent
In similar fashion, if $ P $ contains a circle $ \beta $ such that deleting  $ \beta $
from $ P $ produces a complementary component of type $ (0,4) $,  then we obtain a new
pants decomposition $ P' $ by replacing  $ \beta $ by a circle $ \gamma $ intersecting
$\beta$ transversely in two points and disjoint from the other curves of $ P $. The 
transformation $ P \to P ' $ in this case is called an {\it  associativity move\/} or
A-move.  (In the surface of type $ (0,4) $  containing $ \beta $ and $ \gamma $ these
two curves separate the  four boundary circles in two different ways, and one can view
these  separation patterns as analogous to inserting parentheses via  associativity.)
Note that the inverse of an S-move is again an S-move, and the inverse of an A-move is
again an A-move. 

\medskip
\noindent {\bf Definition.} The {\it pants decomposition complex\/}  $ \P(\Sigma) $ is the two-dimensional cell complex having  vertices the isotopy classes of pants decompositions of $ \Sigma $,  with an edge joining two vertices whenever the corresponding maximal  cut systems differ by a single S-move or A-move, and with faces added to  fill in all cycles of the following five forms: 
{\parindent=27pt
\smallskip
\item{(3A)} Suppose that deleting one circle from a maximal cut system  creates a complementary component of type $ (0,4) $. Then in this  component there are circles $ \beta_1 $, $ \beta_2 $, and $ \beta_3 $,  shown in Figure 2(a), which yield a cycle of three A-moves: $ \beta_1 \to  \beta_2 \to \beta_3 \to \beta_1 $. (No other loops in the given  pants decomposition change.)
\vskip.1truein
\centerline{\epsfbox{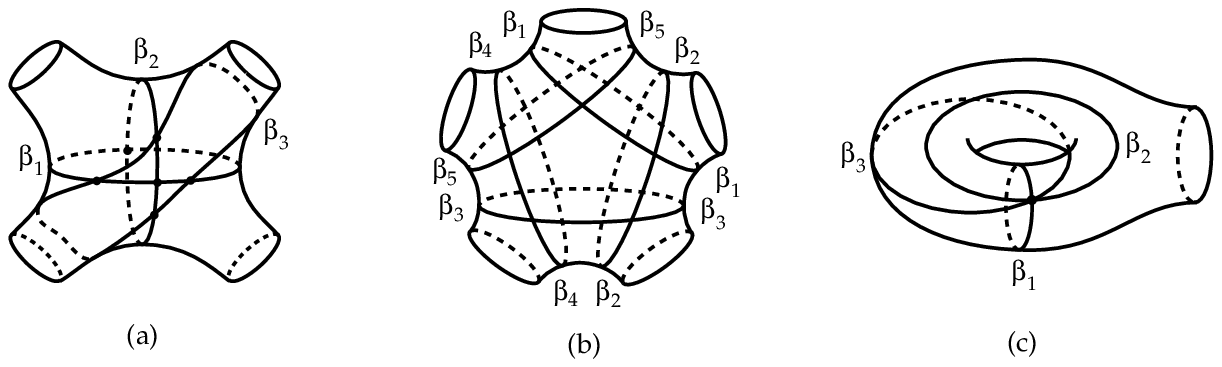}}
\centerline {\bf Figure 2}

\vskip .2cm
\item{(5A)} Suppose that deleting two circles from a maximal cut system creates  a complementary component of type $ (0,5) $. Then in this component  there is a cycle of five A-moves involving the circles $ \beta_i
$  shown in Figure 2(b): $ \{\beta_1, \beta_3\} \to \{\beta_1,\beta_4\} \to 
\{\beta_2,\beta_4\} \to \{\beta_2,\beta_5\} \to \{\beta_3,\beta_5\}  \to
\{\beta_3,\beta_1\} $.

\vskip .3cm
\item{(3S)}Suppose that deleting one circle from a maximal cut system creates a
complementary component of type $ (1,1) $. Then in this  component there are circles $
\beta_1 $, $ \beta_2 $, and $ \beta_3 $, shown in Figure 2(c), which yield a cycle of
three  S-moves: $ \beta_1 \to \beta_2 \to \beta_3 \to \beta_1 $. 

\vskip .3cm
\item{(6AS)}Suppose that deleting two circles from a maximal cut system  creates a complementary component of type $ (1,2) $.  Then in this  component there is a cycle of four A-moves and two S-moves shown in  Figure 3: $\{\alpha_1,\alpha_3\} \to \{\alpha_1,\varepsilon_3\}\to \{\alpha_2,\varepsilon_3\}\to 
\{\alpha_2,\varepsilon_2\}\to \{\alpha_2,\varepsilon_1\}\to \{\alpha_3,\varepsilon_1\}\to \{\alpha_3,\alpha_1\}$. 

\vfill\eject
\vskip .3cm
\centerline{\epsfbox{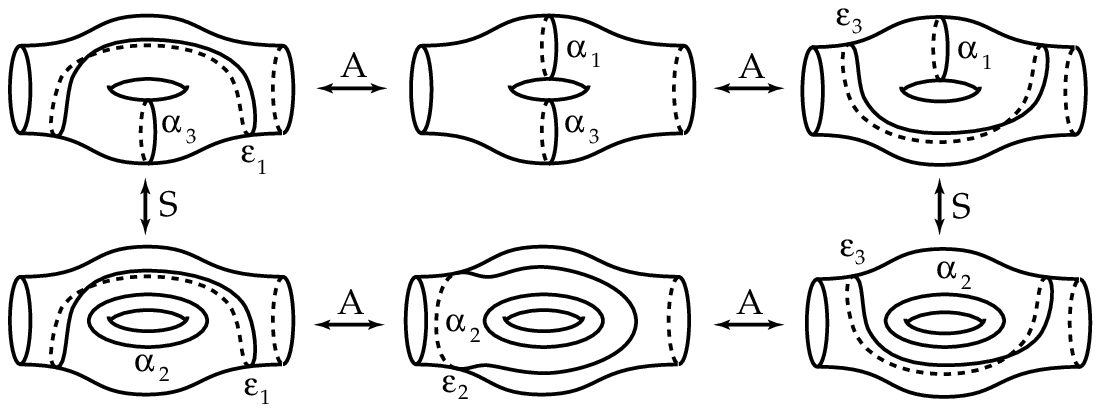}}
\vskip .3cm
\centerline {\bf Figure 3}

\vskip .3cm
\item{(C)} If two moves which are either A-moves or S-moves are  supported in  disjoint
subsurfaces of $ \Sigma $, then they commute,  and their commutator is a cycle of four
moves.  

}

\vskip .3cm
\noindent {\bf Theorem.} {\sl The pants decomposition complex $ \P(\Sigma) $ is simply connected.}
\vskip .2cm
Thus any two sequences of A-moves and S-moves joining  two given pants decompositions can be obtained one from the other by a finite  number of insertions or deletions of the five types of cycles,  together with the trivial operation of inserting or deleted a  move followed by its inverse.  

\medskip
\hfill\smash{\lower 142pt\llap{\epsfbox{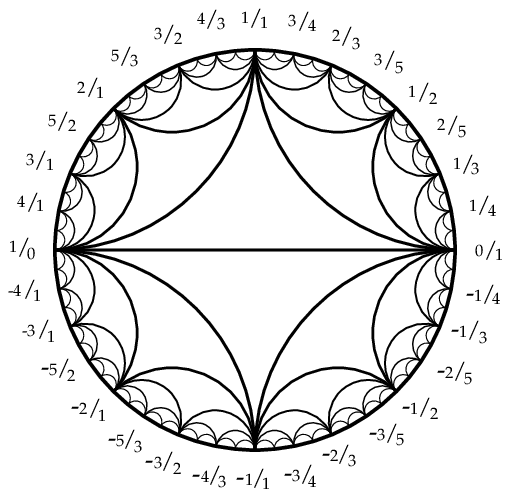}}} 
\vskip-\baselineskip
\hangafter=-12\hangindent=-13pc
\noindent{\bf Example.} If $ \Sigma $ has type (0,4) or  (1,1), the
two cases when a maximum cut system contains just one circle,  then $ \P(\Sigma) $ is
the two-dimensional complex shown in Figure 4,  consisting entirely of triangles since
only the relations 3A or 3S are possible in these two cases. The vertices of  $
\P(\Sigma) $ are labelled by slopes, which classify the  nontrivial isotopy classes of
circles on $ \Sigma $.  This is a familiar fact for the torus, where slopes are defined
via homology. For the $(0,4)$ surface, slopes are defined by lifting curves to the torus
via the standard two-sheeted branched covering of the sphere by the torus, branched over
four points which become the four boundary circles of the $(0,4)$ surface.

\hfill\smash{\raise24pt\llap{\bf Figure 4 \qquad\qquad\ }}

\noindent P{\tiny ROOF}. This uses the same basic  approach as in [HT], which consists of realizing cut systems as  level sets of Morse functions $ f \: \Sigma \to \R$. 

Let $I=[0,1]$.  We consider Morse functions $ f \: (\Sigma, \bdy\Sigma) \to  (I,0) $ whose critical points all lie in the interior of $ \Sigma $. To such a  Morse function we associate a finite graph $ \Gamma(f) $, which is the  quotient space of $ \Sigma $ obtained by collapsing all points in the same component of a level set $ f^{-1}(a) $ to a single  point in $ \Gamma(f) $. If we assume $ f $ is generic, so that all  critical points have distinct critical values, then the vertices of  $ \Gamma(f) $ all have valence 1 or 3 and arise from critical points of  $ f $ or from boundary components of $ \Sigma $. Namely, boundary  components give rise to vertices of valence 1, as do local maxima and  minima of $ f $, while saddles of $ f $ produce vertices of valence 3. See Figure 2 of [HT] for pictures. We can associate to  such a function $ f $ a maximal cut system $ C(f) $, unique up to  isotopy, by either of the following two equivalent procedures: 
\smallskip
\item{(1)} Choose one point in the interior of each edge of $ \Gamma(f) $,  take the circles in $ \Sigma $ which these points correspond to, then  delete those circles which bound disks in $ \Sigma $ or are isotopic to  boundary components, and replace collections of mutually isotopic circles  by a single circle.

\item{(2)} Let $ \Gamma_0(f) $ be the unique smallest subgraph  of $ \Gamma(f) $ which $ \Gamma(f) $ deformation retracts to and  which contains all the vertices corresponding to boundary components  of $ \Sigma $. If $ \Gamma_0(f) $ has vertices of valence 2, regard these  not as vertices but as interior points of edges. In each edge of $ \Gamma_0(f) $ not having a valence 1 vertex as an endpoint,  choose an interior point distinct from the points which were vertices  of valence 2. Then let $ C(f) $ consist of the circles in $ \Sigma $  corresponding to these chosen points of $ \Gamma_0(f) $.

\smallskip\noindent
Every maximal cut system arises as $ C(f) $ for some generic  $ f \: (\Sigma, \bdy\Sigma) \to (I,0) $. To obtain such an $ f $,  one can first define it near the circles of the given cut system  and the circles of  $ \bdy\Sigma $ so that all these circles are noncritical level curves, then extend to a function defined on all of $ \Sigma $,  then perturb this function to be a generic Morse function.

\medskip 
After these preliminaries, we can now show that $ \P(\Sigma) $ is connected. Given two maximal cut systems, realize them as $ C(f_0) $ and $ C(f_1) $.  Connect the generic Morse functions $ f_0 $ and $ f_1 $ by a  one-parameter family $ f_t \: (\Sigma, \bdy\Sigma) \to (I,0) $ with no  critical points near $ \bdy \Sigma $. This is possible since the  space of such functions is convex. After perturbing the family $ f_t $ to be generic, then $ f_t $ is a generic Morse function for each $ t $,  except for two phenomena: birth-death critical points, and {\it crossings\/} interchanging the heights of two consecutive nondegenerate critical  points, as described on p.224 of [HT]. The associated maximal  cut systems $ C(f_t) $ will be independent of $ t $ except for  possible changes caused by these two phenomena. Birth-death  points are local in nature and occur in the interior of an  annulus in $ \Sigma $ bounded by two level curves, hence produce no change in $ C(f_t) $. Crossings can affect $ C(f_t) $  only when both critical points are saddles. Up to level-preserving  diffeomorphism, there are five possible configurations for such a  pair of saddles, shown in Figures 5 and 6 of [HT]. The three  simplest configurations are shown in Figure 5 below, and one can  see that the intermediate level curve dividing the subsurface into  two pairs of pants changes by an A-move as the relative heights of the  two saddles are switched.

\vskip.1cm
\centerline{\epsfbox{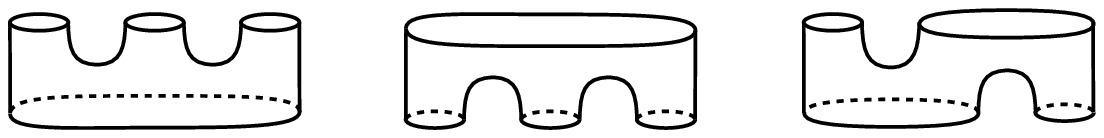}}
\vskip .2cm
\centerline {\bf Figure 5}
\vskip .2cm

The fourth configuration, shown in the left picture of Figure 6 below,  also occurs in
a subsurface of type (0,4).  Here the crossing produces  an interchange of the level
curves $ \alpha_1 $ and $ \alpha_2 $ indicated in the middle picture.  These two curves
intersect in four points, and can be redrawn as in the  right picture. They are related
by a pair of A-moves,  interpolating between them the horizontal circle $ \beta $. (In
terms of  Figure 4, we can connect the slope $ 1 $ and $ -1 $ vertices by an  edgepath
passing through either the slope $ 0 $ or slope $ \infty $  vertices.) 

\centerline{\epsfbox{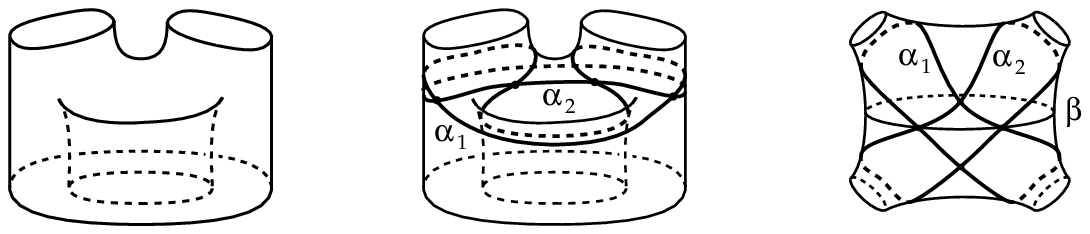}}
\vskip .2cm
\centerline {\bf Figure 6} 
\vskip .3cm

The fifth configuration takes place in a subsurface of  type (1,2), as shown in Figure 7.  Here the two level curves in the 
left-hand figure change to the two in the right-hand figure. This is  precisely the
change from the pair of circles in the middle of the  upper row of Figure 3 to the pair
in the middle of the lower row.  Thus the change is realized by an A-move, an S-move,
and an A-move.  This finishes the proof that $ \P(\Sigma) $ is connected. 

\vskip.2cm
\centerline{\epsfbox{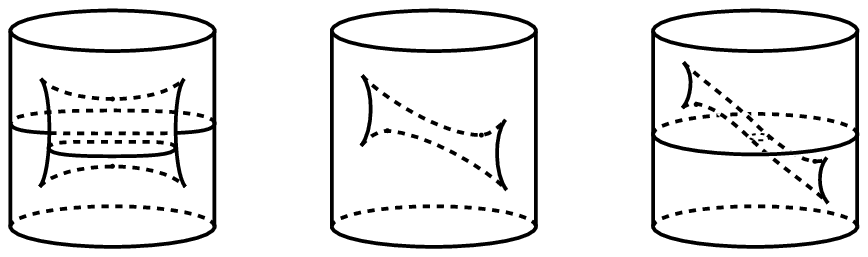}}
\vskip .2cm
\centerline {\bf Figure 7}
\vskip .2cm

Note that the edgepath in $ \P(\Sigma) $ associated to the generic family $ f_t $ is not quite unique. For a crossing as in the  fourth configuration, shown in Figure 6, there were two associated  edgepaths in $ \P(\Sigma) $, which in Figure 4 corresponded to  passing from slope $ 1 $ to slope $ -1 $ through either slope $ 0 $ or slope $ \infty $. These two edgepaths are  homotopic in $ \P(\Sigma) $ using two relations of type 3A.  Similarly, a crossing in the fifth configuration, in Figure 7,  corresponded to an edgepath of three edges, but there are  precisely two choices for this edgepath, the two ways of going halfway around Figure 3, so these two choices are related by a relation of type 6AS. Thus we conclude that the  edgepath associated to a generic family $ f_t $ is unique up to  homotopy in $ \P(\Sigma) $.

A preliminary step to showing $ \P(\Sigma) $ is simply connected is: 

\vskip .3cm
\noindent {\bf Lemma}. {\sl Every edgepath in $ \P(\Sigma) $ is homotopic in  the 1-skeleton of $ \P(\Sigma) $ to an edgepath which is the sequence  of maximal cut systems $ C(f_t) $ associated to a generic one-parameter  family $ f_t $. }

\vskip .2cm
\noindent P{\tiny ROOF}. First we show:
\smallskip
\item{($*$)} If the cut systems $ C(f_0) $ and $ C(f_1) $ are  isotopic, then there is a generic family $ f_t $ joining $ f_0 $ and  $ f_1 $ such that for all $ t $, $ f_t $ has nonsingular level curves in the  isotopy classes of all the circles of $ C(f_0) $ and $ C(f_1) $. 
\smallskip\noindent
This can be shown as follows. Composing $ f_0 $ with an ambient  isotopy of $ \Sigma $ taking the curves in $ C(f_0) $ to the  curves in $ C(f_1) $, we may assume that $ C(f_0) = C(f_1) $. The  normal directions to these curves defined by increasing values  of $ f_0 $ and $ f_1 $ may not agree, but this can easily be  achieved by a deformation of $ f_0 $ near $ C(f_0) $. Then we can  further deform $ f_0 $ so that it agrees with $ f_1 $ near $ C(f_0) =  C(f_1) $ and near $ \bdy \Sigma $, without changing the local  behavior near critical points. Then, keeping the new $ f_0 $  fixed where we have made it agree with $ f_1 $, we can deform it  to coincide with $ f_1 $ everywhere by a generic family $ f_t $. 

To  deduce the lemma from ($*$) it then suffices to realize an  arbitrary A-move or S-move. For A-moves we can just use Figure 5. Similarly, Figure 7 realizes a given S-move sandwiched between  two A-moves, but we can realize the inverses of these A-moves,  so the result follows.  \qed

 \medskip 
Now consider an arbitrary loop in $ \P(\Sigma) $. By the lemma, together with the statement ($*$) in its  proof, this loop is homotopic to a loop of the form $ C(f_t) $  for a loop of generic functions $ f_t $. Since the space of  functions is convex, there is a 2-parameter family $ f_{tu} $  giving a nullhomotopy of the loop $ f_t $. We may assume $ f_{tu} $  is a generic 2-parameter family, so that $ f_{tu_0} $ is a generic 1-parameter family for each $ u_0 $ except for the six types of isolated phenomena  listed on page 230 of [HT]. The proof that $ \P(\Sigma) $ is  simply connected will be achieved by showing that these phenomena  change the associated loop $ C(f_{tu_0}) $ by homotopy in $ \P(\Sigma) $.

The first three of the six involve degenerate critical points and are  uninteresting for our purposes. In each case the change in  generic 1-parameter family is supported in subsurfaces of $ \Sigma $ of  type $ (0,k) $, $ k \leq 3 $, bounded by level curves, so there is no  change in the associated path in $ \P(\Sigma) $.

The last three phenomena, numbered (4), (5), and (6) on page 230 of [HT], involve  only
nondegenerate critical points, which we may  assume are saddles since otherwise the
reasoning in the  preceding paragraph shows that nothing interesting is happening. 
Number (4) is rather trivial: A crossing and its ``inverse" are  cancelled or
introduced. We may choose the segment of the  edgepath in $ \P(\Sigma) $ associated to
the crossing and its  inverse so that it simply backtracks across up to three edges, 
hence the edgepath changes only by homotopy.  Number (5) is the  commutation of two
crossings involving four distinct saddles. This  corresponds to a homotopy of the
associated edgepath across 2-cells  representing the commutation relation C.  Number (6)
arises when  three saddles have the same $ f_{tu} $-value at an isolated  point in the $
(t,u) $-parameter space.  As one circles around this  value, the heights of the saddles
vary through the six possible orders:  123, 132, 312, 321, 231, 213, 123. To finish the
proof it  remains to analyze the various possible configurations for these  three
saddles. The interesting cases not covered by previous arguments are when the three
saddles lie in a connected subsurface bounded by  level curves just above and below the
three saddles. Note that we can  immediately say that all relations among moves, apart
from the commutation relation, are  supported in subsurfaces of types (0,5) and (1,3). 
This is because a subsurface bounded by level curves with three  saddles, hence Euler
characteristic $ -3 $, must have at least  two boundary circles, one below the saddles
and one above, so if the  surface is connected it must have type (0,5) or (1,3). The
analysis  below will show that the (1,3) subsurfaces can be reduced to (1,2) 
subsurfaces. There are sixteen possible configurations of 
three saddles on one level, shown in Figure 8, where the saddles are 
regarded as 1-handles, or rectangles, attached to level curves. The sixteen
configurations are  grouped into eight pairs, the two configurations in each pair
being related by  replacing $ f_{tu}$ by its negative.
\vskip.2cm
\centerline{\epsfbox{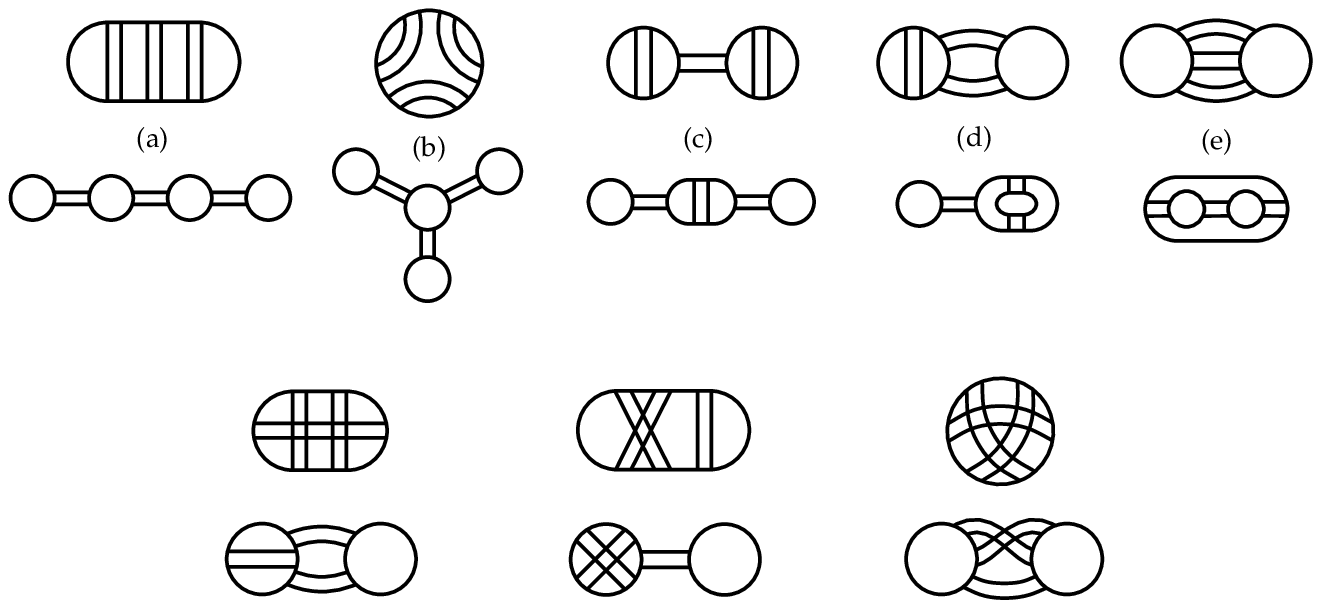}}
\vskip .3cm
\centerline {\bf Figure 8}
\vskip .3cm
\noindent The first five pairs involve a genus zero subsurface and are somewhat 
easier to analyze visually than the other three pairs, which occur in a 
genus one subsurface. We consider each of these five pairs in turn.
\smallskip
\noindent (a) A picture of the subsurface with $ f_{tu} $ as the height  function is shown in Figure 9. 

\centerline{\epsfbox{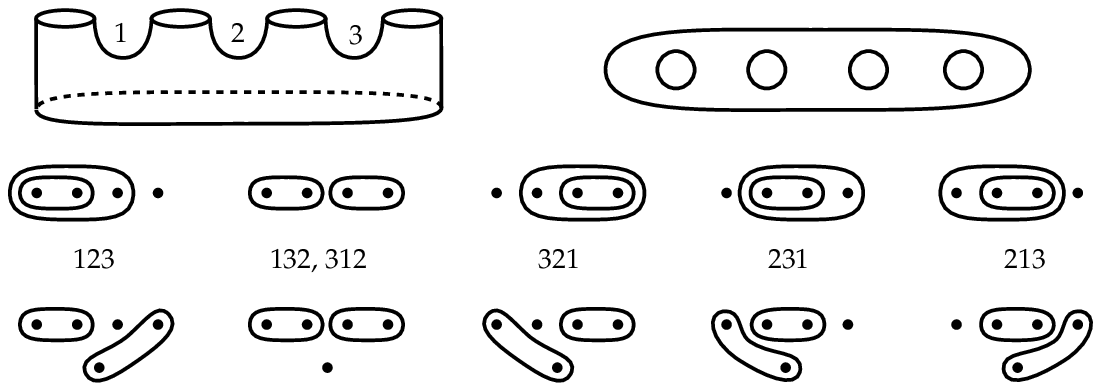}}
\vskip .2cm
\centerline {\bf Figure 9}
\vskip .2cm
\noindent Viewed from above, the surface can be seen as a disk  with  four subdisks
deleted, a (0,5) surface. In the second row of the  figure we show the various
configurations of level curves when the saddles  are perturbed to each of the six
possible orders. For example, the  first diagram shows the order 123, where the saddle 1
is the highest, 2 is the middle, and 3 is the lowest. The two circles shown lie between 
the two adjacent pairs of saddles. The four dots represent four of the  five boundary
circles of the (0,5) surface, the fifth being  regarded as the point at infinity in the
one-point compactification  of the plane. In the third row of the figure this fifth
point is brought  in to a finite point and the level circles are redrawn accordingly. 
The two adjacent orderings 132 and 312 produce the same level curves,  so we have in
reality a cycle of five maximal cut systems. Each is  related to the next (and the first
to the last) by an A-move, and the  whole cycle is the relation 5A.

\smallskip
\noindent (c) We treat this case next since it is very similar to (a). From Figure 10 it is clear that one again has the relation 5A.

\centerline{\epsfbox{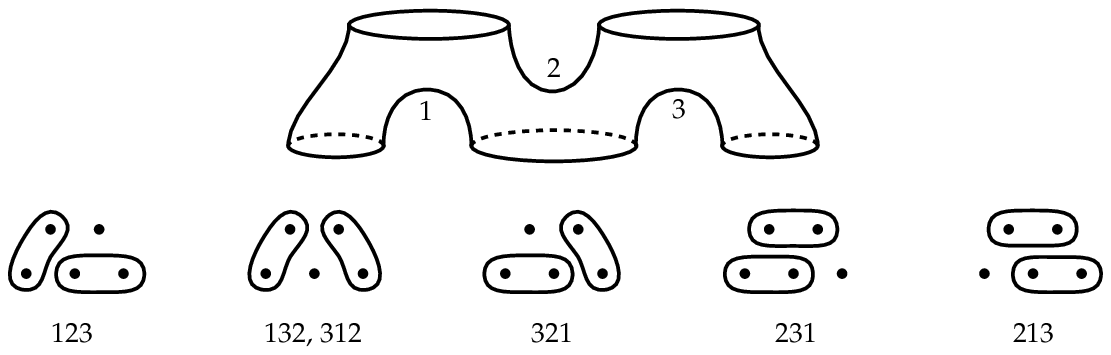}}
\vskip .2cm
\centerline {\bf Figure 10}
\vskip .2cm

\vfill\eject

\noindent (b) Here the 3-fold rotational symmetry makes it unlikely  that one would directly get the relation 5A. The second row of Figure 11  shows the cycle of six cut systems. 

\centerline{\epsfbox{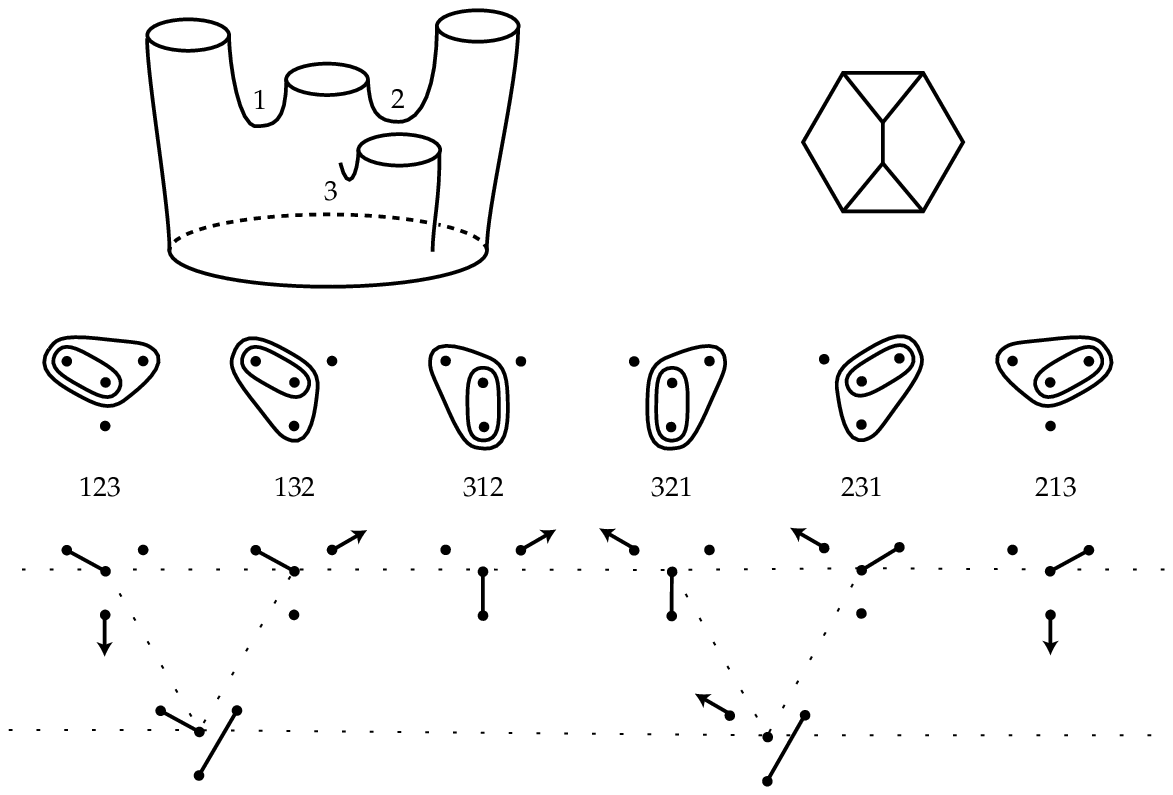}}

\centerline {\bf Figure 11}
\smallskip

\noindent It is convenient to simplify the notation at this point by  representing the two circles in a pants decomposition of the (0,5)  surface by two arcs joining four of the five points representing the  boundary circles. The boundary of a neighborhood of each arc is then a  circle separating two of the five points from the other three. The  third row of the
figure shows the cycle of six cut systems in this  notation, with the fifth point at
infinity and an arc to this point indicated by an arrow from one of the other four
points. Note  that we have a cycle of six A-moves. This can be reduced to  two 3A and
two 5A relations by adjoining the two configurations in the  fourth row.  Schematically,
one subdivides a hexagon into two pentagons  and two triangles by inserting two interior
vertices, as shown. 
\smallskip

\noindent (d) Here the cycle of six cut systems contains two steps  which are not A-moves but resolve into a pair of A-moves. Thus we  have a cycle of eight A-moves, and this decomposes into two 5A  relations, as shown in Figure 12.

\smallskip
\noindent(e) In this case we have the configuration shown in Figure 13, with 
3-fold symmetry.  The cycle of six multicurves has three steps  
which resolve into pairs of A-moves, so we have a cycle of nine A-moves. 
This can be reduced to three 3A relations and four cycles of six A-moves.
After a permutation of the five boundary circles of the (0,5) surface, 
each of these 6-cycles becomes the 6-cycle considered in case (b).

\vskip.2cm
\centerline{\epsfbox{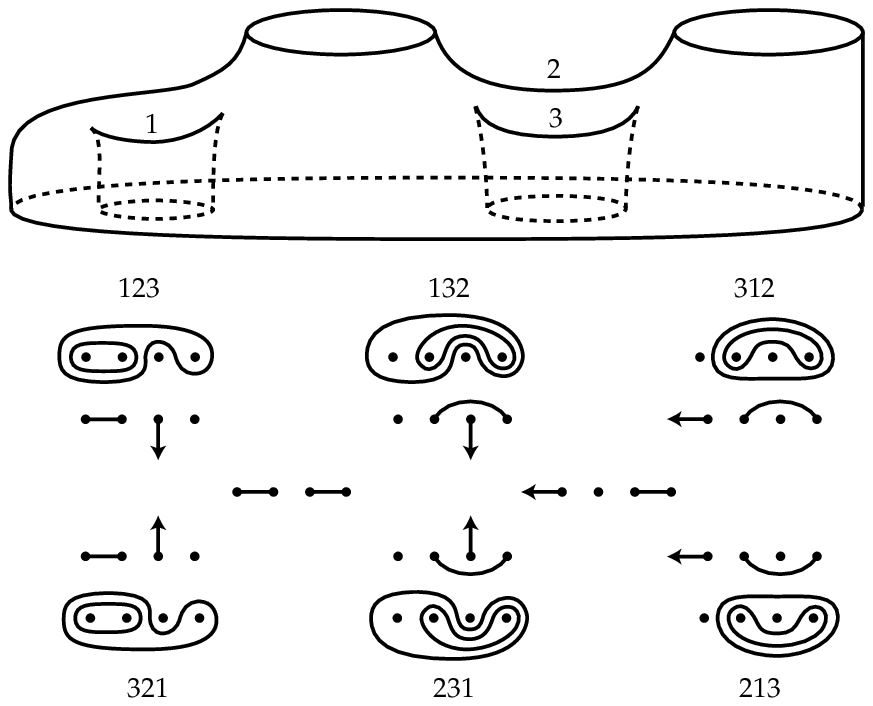}}
\vskip .2cm
\centerline {\bf Figure 12}
\vskip .8cm

\centerline{\epsfbox{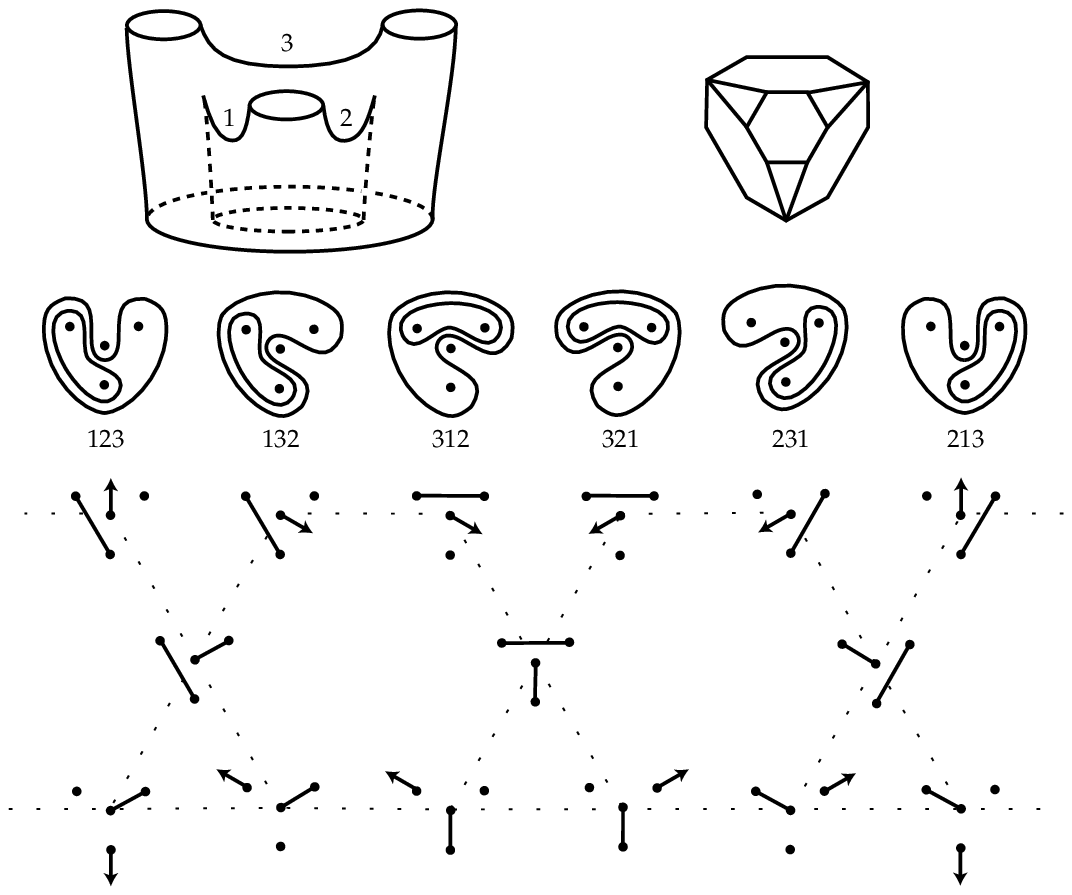}}

\centerline {\bf Figure 13}
\vskip .2cm
\smallskip 
This completes the analysis of the five cases of the phenomenon (5)  involving  genus 0
surfaces. In particular, the theorem is now  proved for surfaces of type $(0,n)$. To
finish the proof it would suffice  to do a similar analysis of the three remaining
configurations of  three saddles in surfaces of type (1,3). However, the cycles of A- 
and S-moves arising from these configurations are somewhat more  complicated than those
in the genus zero configurations, so instead of carrying out this analysis, we shall
make a more general  argument, showing that the relation 3A and 6AS suffice to reduce
the  genus one case to the genus zero case.  So let $ \Sigma $ have  type $ (1,n) $. We
can view the boundary components of $ \Sigma $ as  punctures rather than circles, so $
\Sigma $ is the complement of $ n $  points in a torus $ \widehat\Sigma $. Given an
edgepath loop $ \gamma $  in $ \P(\Sigma) $, its image $ \widehat\gamma $ in $
\P(\widehat\Sigma) $ is  nullhomotopic since the explicit picture of $ \P(\widehat\Sigma)
$  shows it is contractible. Our task is to show the nullhomotopy  of $ \widehat\gamma $
lifts to a nullhomotopy of $ \gamma $. 

The nullhomotopy of $ \widehat\gamma $ gives a map $ \widehat g \: D^2 \to 
\P(\widehat\Sigma)
$. Making $ \widehat g $ transverse to the graph dual to the  1-skeleton of $
\P(\widehat\Sigma) $, the preimage of this dual  graph is a graph $ G $ in $ D^2 $,
intersecting the boundary of $ D^2 $  transversely, as depicted by the solid lines in
the left half of Figure 14.

{\epsfxsize=11cm\centerline{\epsfbox{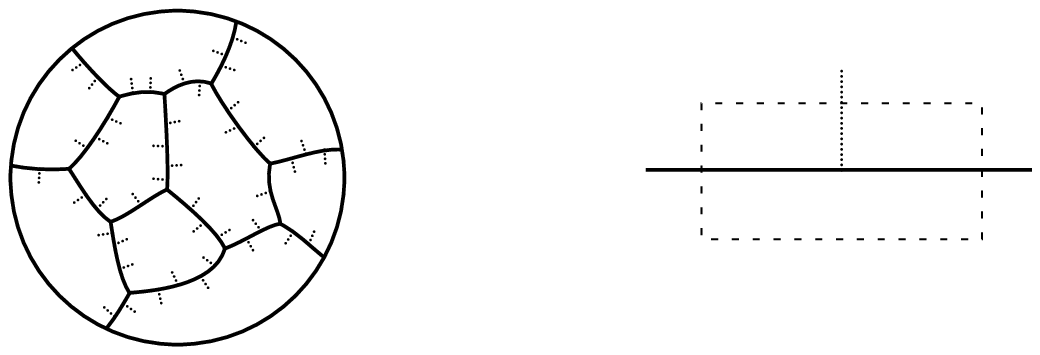}}}
\centerline {\bf Figure 14}
\vskip .2cm
\noindent The vertices of $ G $ in the interior of $ D^2 $ have  valence three,  and are
the preimages of the center points of  triangles of $ \P(\widehat\Sigma) $. Each such
vertex corresponds to  three circles on $ \widehat\Sigma $ having distinct slopes and 
disjoint except for a single point where they all three intersect transversely. Each
edge of $ G $ corresonds to a pair of circles on $ \widehat\Sigma $  of distinct slopes,
intersecting transversely in one point. The  complementary regions of $ G $ correspond
to single circles. 

In a neighborhood $ N $ of $ G $ we can choose all these circles in $ \widehat\Sigma $
to vary continuously with the point in $ N $. We  can also assume these continuously
varying circles have general position  intersections with the $ n $ puncture points, so
that they are disjoint  from the punctures except along arcs, shown dotted in Figure 14,
abutting  interior points of edges of $ G $, where a single circle slides across a 
puncture. Near such a dotted arc we thus have three circles: the circle before  it
slides across the puncture, the circle after it slides across the  puncture, and a third
circle intersecting each of the two circles in one  point transversely.  We can perturb
the first two circles to be disjoint,  so they are essentially two parallel copies of
the same circle with the  puncture between them. A neighborhood of the three circles is
then a surface  of type (1,2). We can identify the three circles in this subsurface with
the  three simplest circles in Figure 3: the upper and lower meridian circles  and the
longitudinal circle. The puncture is one of the two boundary circles of the subsurface.
Adjoining the other circles shown in the figure, we get  various pants decompositions of
the subsurface. Choosing a fixed  pants decomposition of the rest of $ \Sigma $ then
gives a way of  lifting $ \widehat g $ to $ g \: D^2 \to \P(\Sigma) $ in a neighborhood of
the  dotted arc, by superimposing Figure 3 on the right half of Figure~14.  Since the
chosen circles are disjoint from punctures elsewhere along $ G $,  we can then extend
the lift $ g $ over $ G $ by extending the given circles to pants decompositions of $
\Sigma $, using just the fact that any  two pants decompositions of a genus zero surface
can be connected by a  sequence of A-moves. Finally, the lift $ g $ can be extended over
the  complementary regions of $ G $ since the theorem is already proved for  genus zero
surfaces.\qed  
\medskip
\noindent {\bf References}
\vskip .3cm
\cond{[HT]} A. Hatcher and W. Thurston, A presentation for the mapping
class group of a closed orientable surface, {\it Topology} {\bf 19} (1980),
221-237.

\end